\documentstyle[12pt]{article}
\title{THE VON NEUMANN REGULAR RADICAL AND JACOBSON RADICAL  OF
CROSSED PRODUCTS
\thanks {This work is supported by National Science Foundation}}
\author{Shouchuan Zhang \\ Department  of Mathematics,
Nanjing    University, 210008 } 
\date{}
\begin{document}
\newtheorem{Theorem}{\quad Theorem}[section]
\newtheorem{Proposition}[Theorem]{\quad Proposition}
\newtheorem{Definition}[Theorem]{\quad Definition}
\newtheorem{Corollary}[Theorem]{\quad Corollary}
\newtheorem{Lemma}[Theorem]{\quad Lemma}
\newtheorem{Example}[Theorem]{\quad Example}
\maketitle
\begin {abstract}

We construct the $H$-von Neumann regular radical for $H$-module algebras
and show that it is an $H$-radical property.  We obtain that
 the Jacobson radical of  twisted graded
algebra  is a graded ideal.  For twisted $H$-module algebra
$R$, we also show that
$r_{j}(R\# _\sigma  H)= r_{Hj}(R)\# _\sigma H$  and the Jacobson radical
of $R$ is stable,
when $k$ is an algebraically closed field or there exists an algebraic
closure $F$ of $k$ such that
$r_j(R \otimes  F) = r_j(R) \otimes F$,
where $H$ is a finite-dimensional, semisimple, cosemisimple, commutative or
cocommutative Hopf algebra
 over $k$. In particular, we answer two questions
 J.R.Fisher asked.
   \end {abstract}
\noindent
\addtocounter{section}{-1}

 \section{ Introduction }
             Throughout this paper, if we do not specially announce,
             we let $k$ be
 a commutative associative ring with unit,
$H$  an
algebra with unit and comultiplication $\bigtriangleup$
( i.e. $\Delta $ is a $k$-module homomorphism: $H \rightarrow H \otimes H$, and coassociative
law holds ),
  $R$  an algebra
over $k$ ($R$ may be without unit).

We define some necessary concepts as follows:

If   there exist a $k$-module homomorphism $
 \left \{ \begin {array} {ll} H \otimes R & \longrightarrow R \\
 h \otimes r & \mapsto   h \cdot r \end {array} \right. $ such that
     $$h\cdot rs = \sum (h_{1}\cdot r)(h_{2}\cdot s) \hbox { \ \ and \ }
     1_H \cdot r = r $$
     for all $r, s \in R, h\in H,$  then we say that
     $H$ weakly acts on $R,$  where $\Delta (h) = \sum h_1 \otimes h_2.$
  For any ideal $I$ of $R$, set
  $$(I:H):= \{ x \in R \mid h\cdot x \in I \hbox { for all } h \in H  \}. $$
 $I$ is called an $H$-ideal,
if $h\cdot I
\subseteq I$ for any $h \in H$. Let $I_H$ denote the maximal $H$-ideal of $R$
in $I$. It is clear that $I_H= (I:H).$
 $R$ is called  $H$-semiprime, if there are no
non-zero nilpotent $H$-ideals in $R$.  $R$ is called  $H$-prime
if $IJ =0$ implies $I=0$ or $J=0$
for any $H$-ideals $I$ and $J$ of $R$.
An $H$-ideal $I$ is called an   $H$-(semi)prime ideal of $R$ if
$R/I$ is  $H$-(semi)prime.   
   $ \{ a_{n} \}$  is called an
  $H$-$m$-sequence  in  $R$ with beginning $a$
  if there exist $h_n, h_n' \in H, b_n \in R$ such that $a_{1} = a \in R$ and
  $a_{n+1} =
  (h_{n}.a_{n})b_{n}(h_{n}'.a_{n})$   for any natural number $n$.
If for every $H$-$m$-sequence $\{ a_{n} \}$
  with  $a_{1} = a$, there exists a natural number $k$ such that $a_{k}= 0,$
  then $a$ is called $H$-$m$-nilpotent element, and we set
  $$W_{H}(R) = \{ a \in R \mid a \hbox { \ is an \ } H \hbox{-}m
  \hbox {-nilpotent element} \} $$  
 $R$ is called an $H$-module algebra
 if the following
conditions hold:

   (i)  $R$ is a unital left $H$-module \ (i.e. $R$ is a left $H$-module and
   $1_H \cdot a = a$
   for any $a \in R$);

   (ii)  $h\cdot ab = \sum (h_{1}\cdot a)(h_{2}\cdot b)$ for any $a, b \in R$,
      $h\in H$, where
$\Delta (h) = \sum h_{1} \otimes h_{2}$.

$R$ is called a twisted $H$-module algebra if
 the following
conditions are satisfied:

   (i) $H$ weakly acts on  $R$;

(ii) $R$ is a twisted $H$-module, that is,
 there exists a  map
 $\sigma \in Hom_k (H \otimes H, R) $ such that
$h \cdot (k \cdot r) = \sum \sigma (h_1,k_1)(h_2k_2 \cdot r)
\sigma ^{-1}(h_3,k_3)$ for all $h, k \in H$ and $r \in R$.

  It is clear that if $H$ has a  counit and $\sigma $  is trivial, then
  the twisted $H$-module algebra $R$ is an $H$-module algebra.

 A left $R$-module $M$ is called an $R$-$H$-module if $M$ is also a
left unital $H$-module
 with $h  (am)= \sum (h_1 \cdot a)(h_2m)$  for all
$h \in H, a \in R, m \in M$. An $R$-$H$-module $M$ is called an
$R$-$H$- irreducible module if there are no non-trivial $R$-$H$-submodules
in $M$ and $RM \not=0$. An algebra homomorphism $\psi :
 (R, \sigma ) \rightarrow (R', \sigma ' )$ is called a
twisted $H$-homomorphism  if $\psi (\sigma ) = \sigma '$ and
 $\psi (h \cdot a) = h \cdot \psi (a)$ for any $h \in H,
a \in R.$   If $\psi $ is surjective we write
$(R, \sigma ) \stackrel {\psi } {\sim } (R', \sigma ' )$.
If $\psi $ is bijective we write  $(R, \sigma ) \stackrel {\psi } {\cong }
 (R', \sigma ' ).$
Let $r_b, r_j, r_l, r_{bm}, r_k$ and $r_n$  denote the Baer radical, the Jacobson
radical, the locally nilpotent radical,
 the Brown-MacCoy radical, the nil radical
 and von Neumann regular radical of algebras, respectively.   Let
   $I \lhd_H R$ denote that $I$ is an $H$-ideal of $R.$

       J.R.Fisher \cite{F} built up the general theory of $H$-radicals for
$H$-module algebras,
studied the $H$-Jacobson radical
and  obtained
\begin {eqnarray}
  r_{j}(R\#H) \cap R  = r_{Hj}(R)
 \label {e (1)}
 \end {eqnarray}
\noindent for any irreducible Hopf algebra $H$(\cite [Theorem 4]{F}).
 J.R.Fisher \cite{F} asked when is
 \begin {eqnarray}             
 r_{j}(R\#H) =r_{Hj}(R) \#H
\label {e (2)}
\end {eqnarray}                
\noindent and asked whether
\begin {eqnarray}
 r_{j}(R\#H) \subseteq (r_{j}(R):H) \#H
\label {e (3)}                          
\end {eqnarray}

If $H$ is a finite-dimensional semisimple Hopf algebra,
then   relation
 (\ref {e (2)})  holds in  the following  three cases:

(i)  the weak action of $H$ on $R$ is inner;

(ii)  $H$ is commutative cosemisimple;

(iii) $H= (kG)^*$, where $G$ is a finite group.

Parts (i) and (ii) are due to author \cite [Thoerem 3.4 and 3.5] {Z98}.
Part (iii)  is due to   M.Cohen and
 S.Montgomery \cite[Theorem 4.1] {CM}.

The author in \cite [Proposition 3.1] {Z98} obtained that
the relation (\ref {e (1)})  holds  for any Hopf algebra $H$.
 In the same paper, the author
proved that relations (\ref {e (2)}) and (\ref {e (3)}) are equivalent
when $H$ is finite-dimensional.
                        M.Cohen and
 S.Montgomery in  \cite[Theorem 4.4] {CM} proved that  the Jacobson radical of graded
 algebra of  finite group type is graded.
                                  
             In this paper,
             we construct the $H$-von Neumann regular radical for $H$-module algebras
and show that it is an $H$-radical property.  We obtain that
 the Jacobson radical of  a twisted graded
algebra  is a graded ideal.  For a twisted $H$-module algebra
$R$, we also show that
$r_{j}(R\# _\sigma  H)= r_{Hj}(R)\# _\sigma H$   and the Jacobson radical
of $R$ is stable by weak action of $H$,
when $k$ is an algebraically closed field or there exists an algebraic
closure $F$ of $k$ such that
$r_j(R \otimes  F) = r_j(R) \otimes F$,
where $H$ is a finite-dimensional, semisimple, cosemisimple, commutative or
cocommutative Hopf algebra
 over $k$. In particular, we answer two questions J.R.Fisher asked in
 \cite {F}.
 Namely, we give a necessary and sufficient condition for
  $r_{j}(R\#H) =r_{Hj}(R) \#H.$ Meantime, we also give an example to show
  that relation (3) does not hold.

\section {The general theory of $H$-radicals for twisted $H$-module algebras }

In this section we give the general theory of $H$-radicals for twisted
$H$-module algebras

\begin {Definition}\label {1.1}
Let $r$ be a property of $H$-ideals of twisted $H$-module algebras. An $H$-
ideal $I$ of twisted $H$-module algebra $R$ is called   an $r$-$H$-ideal of $R$
if it is of the $r$-property.
A twisted $H$-module algebra $R$  is called   an $r$-twisted $H$-module algebra
if it is $r$-$H$-ideal of itself.
A property $r$ of $H$-ideals of twisted $H$-module algebras   is called an
 $H$-radical property if the following conditions are satisfied:

(R1) Every twisted $H$-homomorphic image of $r$-twisted $H$-module algebra
is an $r$ twisted $H$-module algebra;

(R2) Every twisted $H$-module algebra $R$ has  the maximal $r$-$H$-ideal $r(R)$;

(R3)  $R/r(R)$ has not any non-zero $r$-$H$-ideal.

\end {Definition}
We call $r(R)$ the $H$-radical of $R$.

\begin {Proposition}\label {1.2}
Let $r$ be an ordinary hereditary radical property for rings. An
$H$-ideal $I$ of twisted $H$-module algebra
$R$ is called  an $r_H$-$H$-ideal  of $R$ if $I$ is an $r$-ideal of ring $R$.
Then $r_H$ is an $H$-radical property for twisted $H$-module algebras
\end {Proposition}

{\bf Proof.} (R1). If $(R,\sigma )$ is an $r_H$-twisted $H$-module algebra
 and $(R, \sigma )  \stackrel {f} { \sim }  (R', \sigma ')$, then $r(R') =R'$
 by ring theory. Consequently, $R'$  is an $r_H$-twisted $H$-module algebra.

(R2). For any  twisted $H$-module algebra $R$, $r(R)$  is the maximal
$r$-ideal of $R$ by ring theory. It is clear that
$r(R)_H$ is the maximal $r$-$H$-ideal,  which is an $r_H$-$H$-ideal of $R$.
Consequently, $r_H(R) = r(R)_H$ is the maximal $r_H$-$H$-ideal of $R$.

(R3). If $I/r_H(R)$  is an $r_H$-$H$-ideal of $R/r_H(R)$, then
$I$ is an $r$-ideal of algebra $R$ by ring theory. Consequently,
$I \subseteq r(R)$ and $I \subseteq r_H(R).$
$\Box$

\begin {Proposition}\label {1.3}
$r_{Hb}$ is an $H$-radical property
\end {Proposition}

{\bf Proof.} (R1). Let $(R,\sigma )$ is an $r_{Hb}$-twisted $H$-module algebra
 and $(R, \sigma )  \stackrel {f} { \sim }  (R' \sigma ')$.
    For any $x' \in R'$
and  any  $H$-$m$-sequence $\{  a_n' \}$    in $R'$
 with $ a_1' =  x' $,  there exist
 $ b_n' \in  R'$ and $h_n, h_n' \in H$ such that
 $  a_{n+1}' =
 (h_n \cdot  a_n')  b_n' (h_n' \cdot  a_n')$
                           {~}{~}  for any natural number $n$.
Let   $a_1, b_i \in R $  such that $f(a_1) =x'$ and   $f(b_i)=b_i'$
for $i = 1, 2, \cdots .$  Set
 $  a_{n+1} = (h_n \cdot  a_n)  b_n (h_n' \cdot  a_n)$ {~}{~}
                             for any natural number $n$.
 Since $ \{ a_n \}$ is an $H$-$m$-sequence in $R$,
 there exists a natural number $k$ such that $a_k =0$. It is
 clear that $f( a_n)= a_n'$
  for any natural number $n$ by induction. Thus $  a_k' = 0,$
  which implies that $x' $  is an $H$-$m$-nilpotent element.
Consequently, $R'$ is an $r_{Hb}$-twisted $H$-module algebra.

(R2). By \cite [Theorem 1.5] {Z98},
$r_{Hb}(R)= W_H(R) = \{ a \mid a$ is an $H$-$m$-nilpotent element in $R \}$.
Thus $r_{Hb}(R)$ is the maximal $r_{Hb}$-$H$-ideal of $R$.

(R3).  It immediately follows from  \cite [Proposition 1.4] {Z98}. $\Box$

\section {The relations among  radical of $R$  , radical of
  $R \# _\sigma H$, and $H$-radical of $R$ }

In this section we give
the relation among the Jacobson radical $r_j (R)$ of $R$  ,the Jacobson
 radical $r_j(R \# _\sigma H)$ of
  $R \# _\sigma H$, and $H$-Jacobson radical $r_{Hj}(R)$ of $R$.

In this section, let $k$ be a field, $R$ an algebra with unit,
$H$ a Hopf algebra over $k$ and $R \# _\sigma H$ an algebra with unit.
Let $r$ be a hereditary radical property for rings which satisfies
$$r(M_{n \times n} (R)) = M_{n \times n} (r(R))$$ for any twisted
$H$-module algebra $R$.

 Example.  $r_j$, $r_{bm}$ and $r_n$ satisfy the above conditions
by \cite {Sz}. Using \cite [Lemma 2.1 (2)]{Z98},we can easily prove that
$r_b$ and $r_l$   also satisfy the above conditions.

\begin {Definition} \label {2.1}
$\bar r_H (R) := r(R \# _\sigma H) \cap R$  and
$r_H (R) := (r(R):H)$
\end {Definition}

     If $H$ is a finite-dimensional   Hopf algebra and $M = R\#_{\sigma} H$,
     then $M$ is a free right $R$-module  with finite rank by
\cite [Proposition 7.2.11]{Mo}
and $End (M_R) \cong (R \#_{\sigma} H)\#H^*$ by \cite [Corollary 9.4.17]{Mo}.
It follows from  part (a) in the proof of \cite [Theorem 7.2] {MS} that
 there exists  unique bijective  map
 $$\Phi: {\cal I}(R) \longrightarrow {\cal I}(R')$$
such that \ \ \  $\Phi (I)M=MI,$
where $R' = (R \#_\sigma H)\# H^*$  and

 ${ \cal I}(R) = \{ I \mid I \hbox { is an ideal of } R \}$.

\begin {Lemma}\label {2.2}  If $H$ is a finite-dimensional Hopf algebra, then

  $ \Phi (r(R)) = r((R \# _\sigma H) \# H^*).$

\end {Lemma}
{\bf Proof.}  It is similar to the proof of \cite
[lemma 2.1 (5)] {Z98}. $\Box$

\begin {Proposition}\label {2.3}
$\bar r_H (R) \# _\sigma H \subseteq r_{H^*}(R \# _\sigma H)
                           \subseteq r(R \# _\sigma H). $

\end {Proposition}
{\bf Proof.}                          
 We see that
\begin {eqnarray*}
   \bar r_{H}(R) \#_{\sigma} H &=&    (\bar r_{H}(R) \#_{\sigma} 1)(1 \#_{\sigma} H)\\
    &\subseteq &  r(R \#_{\sigma} H) ( 1 \#_\sigma H) \\
    &\subseteq & r(R \#_{\sigma} H). { \ \ \ \  }
\end {eqnarray*}
Thus      $\bar r_H (R) \# _\sigma H \subseteq r_{H^*}(R \# _\sigma H)$
since     $\bar r_H(R) \# _\sigma H $  is an $H^*$-ideal of
$R \# _\sigma H $.   $\Box$

\begin {Proposition}\label {2.4}
If $H$ is a finite-dimensional
      Hopf algebra,  then

(1) $r_H (R) \# _\sigma H =  \bar r_{H^*}(R \# _\sigma H)$;

Furthermore, if               $\bar r_H \le r_H $, then

(2)    $\bar r_H = r_H$ and
$r_H (R) \# _\sigma H \subseteq r(R \# _\sigma H)$;

(3)    $R\# _\sigma H$  is $r$-semisimple
for any $r_H$-semisimple $R$ iff 
   $$r(R \#_\sigma H) = r_H(R) \# _\sigma H.$$

\end {Proposition}
{\bf Proof.} Let $A = R \# _\sigma H.$

(1)         We see that
\begin{eqnarray*}
(r_{H}(R)\#_{\sigma} H)\#H^* &=& \Phi(r_{H}(R))\\
 &=& (\Phi(r(R) )\cap A) \#H^*  \hbox { \ \  by\   \cite [Theorem 7.2]{MS} }\\
&=&(r(A\#H^*)\cap A) \# H^* \hbox {  \  \ \ by  Lemma \ref {2.2}  } \\
&=& \bar r_{H^*}(A)\#H^* \hbox { \ by \ Definition \ref{2.1} }.
\end {eqnarray*}
Thus $ \bar r_{H^*}(A)= r_{H}(R)\#_{\sigma} H.$

(2) We see that
\begin {eqnarray*}
\bar  r_{H}(R) &=& r(A) \cap R \\
& \supseteq & \bar r_{H^*}(A) \cap R \hbox  { \ by assumption }\\
&=&  r_{H}(R)  \hbox { \ by \ part  \  (1) }.
\end {eqnarray*}
Thus  $\bar r_{H}(R)= r_{H}(R)$ by assumption.

(3) Sufficiency is obvious. Now we show the necessity.
            Since $$r((R \# _\sigma H)/(r_H (R) \#_\sigma H))
            \cong r(R/r_H(R) \#_{\sigma '}H) =0 $$
            we have $r (R \#_\sigma H) \subseteq r_H(R) \#_\sigma H.$
Considering part (2), we have

   $$r(R \#_\sigma H) = r_H(R) \# _\sigma H. { \ \ \ } \Box  $$

\begin {Corollary}\label {2.5}
Let $r$ denote $r_b, r_l,  r_j, r_{bm}$ and $r_n.$ Then

(1) $\bar r_H \le r_H;$

Furthermore, if $H$ is a finite-dimensional Hopf algebra, then

(2) $\bar r_H = r_H $;

(3) $R \# _\sigma H $ is $r$-
semisimple for any $r_{H}$-semisimple $R$  iff $r(R \# _\sigma H)=
r_{H}(R)\#_\sigma H$;

(4)  $R \# _\sigma H $ is $r_j$-
semisimple for any $r_{Hj}$-semisimple $R$  iff $r_j(R \# _\sigma H)=
r_{Hj}(R)\#_\sigma H.$

\end {Corollary}
{\bf Proof.}  (1). When $r= r_b$  or $r = r_j$, it was proved  in \cite
[Proposition 2.3 (1) and 3.2 (1)] {Z98}. The others can similarly be proved.

(2). It follows from Proposition \ref {2.4} (2).

(3) and (4)  follow from part (1) and Proposition \ref {2.4} (3).
$\Box$
\begin {Proposition}\label {2.6}
 If $H = kG$ or  the weak action of $H$ on $R$ is inner,
then

(1). $r_H (R) = r(R)$;

(2)  If, in  addition,   $H$ is a finite-dimensional
 Hopf algebra and $\bar r _H \le r _H$, then
 $r_H(R) = \bar r_H(R) = r(R).$

\end {Proposition}

{\bf Proof.}
(1)  It is trivial.

(2) It immediately follows from part (1) and Proposition \ref {2.1} (1)
(2).  $\Box$

\begin {Theorem}\label {2.8}
Let $G$ be a finite group and $\mid G \mid ^{-1} \in k$.
If $H= kG$ or $H= (kG)^*$, then

(1)  $r_j (R) = r_{Hj}(R)= r_{jH}(R)$;

  (2)   $r_{j}(R\# _\sigma  H)= r_{Hj}(R)\# _\sigma H.$
                  \end {Theorem}
{\bf Proof.}
(1) Let $H =kG$.  We  can easly check
$r_j(R) = r_{jH}(R)$  using the method similar to  the proof of
 \cite [Proposition 4.6]  {Z97}. By  \cite [Proposition 3.3 (2) ] {Z98},
 $r_{Hj}(R) =r_{jH}(R)$. Now,
we only need to show that
  $$r_j (R) = r_{H^*j}(R). $$ 
 We see that
\begin {eqnarray*}
  r_{j}((R\# _\sigma  H^*) \# H) &=&
   r_{H^*j}((R\# _\sigma H^*) \# H)   \hbox { \ \ by \cite [Theorem 4.4 (3)]
   {CM} } \\
   &=&  r_{Hj}(R\# _\sigma  H^*) \# H  \hbox { \ \ by \cite [ Proposition
   3.3 (1)] {Z98}}  \\
   &=& (r_{H^*j}(R)\# _\sigma H^* ) \# H  \hbox { \ \ by \cite [ Proposition
   3.3 (1)] {Z98}}.
\end {eqnarray*}
On the one hand, by       \cite [ Lemma 2.2 (8)] {Z98},
$\Phi (r_j (R)) =  r_{j}((R\# _\sigma  H^*) \# H).$
On the other hand, we have that
$ \Phi (r_{H^*j} (R)) =  ( r_{H^*j}(R )\# _\sigma H^*) \# H $  \ \
 by \cite [Lemma 2.2 (2)]    {Z98}.
 Consequently, $r_j (R) = r_{H^*j}(R)$.

(2) It immediately follows from part (1) and \cite [Proposition 3.3 (1) (2)]
{Z98}.
  $\Box$

\begin {Corollary}\label {2.9}
Let $H$ be a semisimple and cosemisimple Hopf algebra over algebraically
closed field $k$. If $H$ is commutative or cocommutative, then
  $$r_j (R) = r_{Hj}(R)= r_{jH}(R)  \hbox { \ \ \ \ and  \ \ \ }
  r_{j}(R\# _\sigma  H)= r_{Hj}(R)\# _\sigma H.$$
                  \end {Corollary}
{\bf Proof.} It immediately follows from Theorem \ref {2.8} and
\cite [Lemma 8.0.1 (c)] {Sw}. $\Box$

    We give an example to show that conditions in Corollary \ref {2.9} can not
be omitted.

\begin {Example}\label {2.10} (see, \cite [Example P20] {Z98})
Let $k$ be a field of characteristic $p > 0 $, $R = k[x]/(x^p)$.
We can define a derivation on $R$ by sending $x$ to $x +1$.
Set $H=u(kd)$, the restricted enveloping algebra, and $A = R\#H.$
Then

(1)
$r_b (A \# H^*) \not= r_{H^*b}(A) \# H^*$;

(2) $r_j (A \# H^*) \not= r_{H^*j}(A) \# H^*$;
               
(3) $r_j (A \# H^*) \not\subseteq  r_{jH^*}(A) \# H^*$.

\end {Example}
{\bf Proof.}        (1).
By \cite [Example P20]  {Z98}, we have $r_b(R) \not=0$ and $r_{bH}(R)=0.$
Since $\Phi (r_b(R) ) = r_b(A \#H^*)\not=0$  and $\Phi (r_{bH}(R))=
r_{bH^*}(A) \#H^*=0,$
we have that  part (1) holds.

(3). We see that $r_j(A \# H^*) = \Phi (r_j(R))$  and
$r_{Hj}(A) \#H^*= \Phi (r_{Hj}(R)).$  Since $R$ is commutative,
$r_j(R)= r_b(R).$  Thus $r_{Hj}(R)= r_{jH}(R)= r_{bH}(R)=0$ and
$r_j(R) = r_b(R) \not=0$, which implies
 $r_j(A\#H^*) \not\subseteq r_{jH^*}(A)\#H^*$.

(2). It follows from part (3).
 $\Box$

This example also answer the question J.R.Fisher asked in \cite {F} :
$$\hbox {Is  \ \ \ } r_j (R \# H) \subseteq  r_{jH}(R) \# H \hbox { \ \ \  ?} $$

If $F$ is an extension field of $k$, we write  $R^F$  for $R \otimes _k F$
(see, \cite [P49 ]{MS}) .

\begin {Lemma}\label {2.11}  If $F$ is an extension field  of $k$, then
                                              
(1)  $H$ is a semisimple  Hopf algebra over $k$  iff $H^F$ is a semisimple
Hopf algebra over $F$;

(2)  Furthermore, if $H$ is a finite-dimensional Hopf algebra, then
  $H$ is a cosemisimple  Hopf algebra over $k$  iff $H^F$ is a cosemisimple
Hopf algebra over $F$.

\end {Lemma}
{\bf Proof.}
(1) It is clear that
$\int _H^l \otimes F = \int _{H^F}^l.$   Thus
 $H$ is a semisimple  Hopf algebra over $k$  iff $H^F$ is a semisimple
Hopf algebra over $F$.

(2)  $(H \otimes F)^* = H^* \otimes F $ since
$H^* \otimes F \subseteq  (H \otimes F)^*$  and $dim_F (H \otimes F ) =
dim_F (H^* \otimes F)= dim _k H$. Thus we can obtain part (2) by  Part (1). $\Box$

By the way, if $H$ is a semisimple Hopf algebra, then $H$ is a seperable algebra by Lemma
\ref {2.11} (see, \cite [P284] {P}).

\begin {Proposition}\label {2.12}
Let $F$ be an algebraic closure of $k$, $R$ an algebra over $k$ and
$$r(R \otimes _k F) = r(R) \otimes _k F,$$
 If $H$ is a finite-dimensional  Hopf algebra  with cocommutative
 coradical over $k$ ,
then
     $$r(R) ^{dim H} \subseteq r_H(R).$$

\end {Proposition}
{\bf Proof.}
It is clear that $H^F$  is   a finite-dimensional  Hopf algebra
 over $F$ and $dim H = dim  H^F =n$. Let $H^F_0$  be the coradical
 of  $H^F$,  $H^F_1 = H^F_0 \wedge H^F_0,
H^F_{i+1} = H^F_0 \wedge H^F_i$ for $i = 1, \cdots, n-1$.
Notice $H^F_0 \subseteq H_0 \otimes F $. Thus $H^F _0 $  is cocommutative.
It is clear that $H^F_0=kG$  by
\cite [Lemma 8.0.1 (c)]{Sw} and $H^F= \cup H^F_i$. It is easy to show that
if $k > i $ then  $$H^F_i \cdot (r(R ^F))^k\subseteq r(R ^F)$$
by induction for $i$. Thus
$$H^F \cdot (r(R ^F))^{ dim H} \subseteq r(R^F)$$
which implies that $(r(R^F))^{dim H }\subseteq r(R^F)_{H^F}$.
By assumption, we have that
$(r(R) \otimes F)^{dim H }\subseteq (r(R) \otimes F)_{H^F}$.
It is clear that  $(I \otimes F)_{H^F} = I_H \otimes F$  for any ideal $I$
of $R$. Consequently,
$(r(R))^{dim H }\subseteq r(R)_H$.   $\Box$

\begin {Theorem}\label {2.13}
Let $H$ be a  semisimple, cosemisimple,
commutative or cocommutative Hopf algebra
 over $k$.
If there exists an algebraic closure $F$ of $k$ such that
$$r_j(R \otimes  F) = r_j(R) \otimes F \hbox { \ \ and \ \ }
r_j((R \# _\sigma H) \otimes  F) = r_j(R \# _\sigma H) \otimes F$$
then
             
(1) $r_j (R) = r_{Hj}(R)= r_{jH}(R); $

 (2)
  $r_{j}(R\# _\sigma  H)= r_{Hj}(R)\# _\sigma H.$
\end {Theorem}
{\bf Proof.}
 (1). By Lemma \ref {2.11},
$H^F$ is semisimple and cosemisimple. Considering Corollary \ref {2.9},
we have that                                               
   $r_j (R^F) = r_{H^Fj}(R^F)= r_{jH^F}(R^F).$
On the one hand, by assumption,
$r_j(R^F) = r_j(R) \otimes F$. On the other hand,
$r_{jH^F} (R^F) = (r_j (R ) \otimes F)_{H^F} = r_{jH}(R) \otimes F$.
Thus $r_j(R) = r_{jH}(R)$.

(2).  It immediately  follows  from part (1). $\Box$

Considering Theorem \ref {2.13} and \cite [Theorem 7.2.13] {P}, we have
                           
\begin {Corollary}\label {2.14}
Let $H$ be a semisimple,  cosemisimple, commutative or cocommutative Hopf
 algebra over $k$. If there exists an algebraic closure $F$ of $k$ such that
$F/k$ is seperable and algebraic, then

  (1) $r_j (R) = r_{Hj}(R)= r_{jH}(R)$;

  (2)   $r_{j}(R\# _\sigma  H)= r_{Hj}(R)\# _\sigma H.$
                  \end {Corollary}

\begin {Lemma}\label {2.15}
(Szasz \cite {Sz})   $$r_j (R) = r_k(R)$$ holds in the following three cases:

(1) Every element in $R$ is  algebraic over $k$ (\cite [Proposition 31.2] {Sz});

(2) The cardinality of $k$  is strictly greater  than the dimension of $R$ and
$k$ is infinite (\cite [Theorem 31.4] {Sz});

(3) $k$ is uncountable and $R$ is  finitely generated (\cite [Proposition 31.5] {Sz}).

\end {Lemma}

\begin {Proposition}\label {2.16}
Let $F$ be an extension of  $k$.
Then
$r(R) \otimes  F \subseteq  r(R \otimes  F),$
where $r$ denote $r_b, r_k, r_l, r_n  .$

\end {Proposition}
{\bf Proof.}
When $r =r_n$, for any $x \otimes a  \in r_n (R) \otimes F$ with $a \not=0$,
there exists $y \in R$ such that
$x =xyx$. Thus $x \otimes a = (x \otimes a)(y \otimes a^{-1}) (x \otimes a)$,
which implies $r_n(R) \otimes F \subseteq r_n(R \otimes F)$.

Similarly, we can obtain the others. $\Box$

\begin {Corollary}\label {2.17}
Let $H$ be a semisimple, cosemisimple, commutative or cocommutative Hopf algebra.
If there exists an algebraic closure $F$ of $k$ such that $F/k$ is a pure
transcendental extension and one of the following three conditions holds:

(i) Every element in $R\#_\sigma H$ is  algebraic over $k$;

(ii) The cardinality of $k$  is strictly greater  than the dimension of $R$ and
$k$ is infinite;

(iii) $k$ is uncountable and $R$ is  finitely generated;

then

(1) $r_j (R) = r_{Hj}(R)= r_{jH}(R)$;

(2)    $r_{j}(R\# _\sigma  H)= r_{Hj}(R)\# _\sigma H$;

(3) $r_j(R)= r_k(R) $  and $r_j(R \#_\sigma H) = r_k(R\#_\sigma H).$
                  \end {Corollary}
{\bf Proof.}
First, we have that part (3) holds by Lemma \ref {2.15}.
We next see that
\begin {eqnarray*}
r_j(R \otimes F) &\subseteq & r_j(R) \otimes F  \hbox { \ \  \cite [Theorem
7.3.4] {P}} \\
&=&  r_k(R) \otimes F  \hbox { \ \  part (3)} \\
&\subseteq & r_k(R \otimes F ) \hbox { \ \ proposition \ref {2.16} } \\
 &\subseteq & r_j(R \otimes F).
\end {eqnarray*}            
Thus $r_j(R \otimes F) = r_j (R) \otimes F.$
Similarly, we can show that                         
 $r_j((R \#_\sigma H) \otimes F) = r_j (R\#_\sigma H) \otimes F.$

Finally, using Theorem \ref {2.13}, we complete the proof. $\Box$

\section {The $H$-Von Neumann regular radical}

In this section, we construct the $H$-von Neumann regular radical
for $H$-module algebras
and show that it is an $H$-radical property.

\begin {Definition} \label {3.1}
Let $a \in R$. If $a\in (H \cdot a) R (H\cdot a)$, then $a$ is called an
$H$-von Neumann regular element, or an $H$-regular element in short.
If every element of $R$ is an $H$-regular, then $R$ is called an $H$-regular
module algebra, written as  $r_{Hn}$-$H$-module algebra.  $I$ is an $H$-ideal
 of $R$ and every element in $I$ is $H$-regular, then $I$ is called an $H$-
 regular ideal.
\end {Definition}

\begin {Lemma}\label {3.2} If $I$ is an $H$-ideal of $R$ and $a \in I$,
then $a$ is $H$-regular in $I$ iff $a$ is $H$-regular in $H$.
\end {Lemma}
{\bf Proof.} The necessity is clear.

Sufficiency: If $a \in (H \cdot a)R (H\cdot a),$ then there exist
$h_i, h_i' \in H, b_i \in R,$  such that
$$a = \sum (h_i \cdot a) b_i (h_i' \cdot a).$$
We see that
\begin {eqnarray*}
a &=& \sum _{i, j} [ h_i \cdot (( h_j \cdot a) b_j (h_j' \cdot a))]b_i
(h_i' \cdot a) \\
&=&      \sum _{i, j} [ ((h_i)_1 \cdot ( h_j \cdot a)) ((h_i)_2 \cdot b_j)
      ((h_i)_3 \cdot (h_j' \cdot a))]b_i
(h_i' \cdot a)    \\
&\in & (H\cdot a) I (H\cdot a).
\end {eqnarray*}
Thus $a$ is an $H$-regular in $I$. $\Box$
     
\begin {Lemma}\label {3.3}
If  $x - \sum _i (h_i \cdot x) b_i (h_i' \cdot x)$ is
$H$-regular, then
$x$ is $H$-regular, where $x, b_i \in R, h_i , h_i' \in H.$
\end {Lemma}
{\bf Proof.}
Since   $x - \sum _i (h_i \cdot x) b_i (h_i' \cdot x)$    is $H$-regular,
there exist $g_i, g_i' \in H, c_i \in R $ such that
$$x - \sum _i (h_i \cdot x) b_i (h_i' \cdot x)  =
\sum _j ( g_j \cdot (x - \sum _i (h_i \cdot x) b_i (h_i' \cdot x)))
c_j ( g_j' \cdot ( x - \sum _i (h_i \cdot x) b_i (h_i' \cdot x))). $$
Consequently, $x \in (H\cdot x)R (H \cdot x).$
$\Box$

\begin {Definition}\label {3.4}
$$r_{Hn}(R):=  \{ a \in R \mid \hbox  { \ the \ }   H\hbox {-ideal } (a)
\hbox { generated by } a \hbox { is } H \hbox {-regular } \}.$$
\end {Definition}

\begin {Theorem}\label {3.5}
$r_{Hn}(R)$ is an $H$-ideal of $R$

\end {Theorem}
{\bf Proof.} We first show that $R r_{Hn}(R) \subseteq r_{Hn}(R).$
 For any $a \in r_{Hn}(R), x \in R,$ we have that
 $(xa) $  is $H$-regular since $(xa) \subseteq (a).$
 We next show that $a-b \in r_{Hn}(R)$  for any $a , b \in r_{Hn}(R).$
 For any $x \in (a-b),$ since $(a-b) \subseteq (a) + (b)$, we have that
 $x = u -v$ and $u \in (a), v\in (b).$
 Say $u = \sum _i (h_i \cdot u) c_i (h_i' \cdot u)$  and $h_i, h_i' \in H,
  c_i \in R.$
We see that
\begin {eqnarray*}
&x& - \sum _i (h_i \cdot x)c_i (h_i' \cdot x)\\
 &=& (u-v) - \sum _i  (h_i \cdot (u-v))c_i (h_i' \cdot (u-v))   \\
&=& -v - \sum _i [ -  (h_i \cdot u)c_i (h_i' \cdot v) -   (h_i \cdot v)c_i (h_i' \cdot u)
+  (h_i \cdot v)c_i (h_i' \cdot v)] \\
& \in & (v) .
\end {eqnarray*}
Thus  $x - \sum _i (h_i \cdot x)c_i (h_i' \cdot x) $ is $H$-regular
and $x$ is $H$-regular by Lemma \ref {3.3}. Therefore
$a-b \in r_{Hn}(R).$ Obviously, $r_{Hn}(R)$ is $H$-stable. Consequently,
$r_{Hn}(R)$ is an $H$-ideal of $R$. $\Box$

\begin {Theorem}\label {3.5}
$r_{Hn}(R/r_{Hn}(R)) = 0.$

\end {Theorem}
{\bf Proof.}   Let $\bar R = r/r_{Hn}(R)$ and  $\bar b = b + r_{Hn}(R) \in r_{Hn}(R/r_{Hn}(R)).$
It is sufficient to show that $b \in r_{Hn}(R).$  For any $a \in (b)$,
it is clear that  $\bar a \in (\bar b)$.  Thus there exist
$h_i, h_i' \in H, \bar c_i \in \bar R$ such that
$$\bar a = \sum _i (h_i \cdot \bar a) \bar c_i (h_i' \cdot \bar a)
= \sum _i \overline { (h_i \cdot  a)  c_i (h_i' \cdot  a)}.$$
Thus  $ a-  \sum _i  (h_i \cdot  a)  c_i (h_i' \cdot  a) \in r_{Hn}(R),$
which implies that $ a $  is $H$-regular. Consequently,   $b \in r_{Hn}(R).$
Namely, $\bar b =0$ and $r_{Hn}(R)=0.$  $\Box$

\begin {Corollary}\label {3.7.1}
$r_{Hn}$ is an $H$-radical property for $H$-module algebras and $r_{nH}
\le r_{Hn}$.
\end {Corollary}
{\bf Proof.} 
(R1). If $R \stackrel {f} {\sim }  R'$  and $R$  is an $r_{Hn}$
-$H$-module algebra then, for any $f(a) \in R'$,
$f (a) \in (H \cdot f(a)) R' (H \cdot f(a)).$  Thus $R'$ is also an $r_{Hn}$-
$H$-module algebra.

(R2). If $I$ is an $r_{Hn}$-$H$-ideal of $R$  and $ r_{Hn}(R) \subseteq I$
then, for any $a\in I,$  $(a)$  is $H$-regular
since $(a) \subseteq I.$  Thus $I \subseteq r_{Hn}(R).$

  (R3). It follows from Theorem \ref {3.6}.

  Consequently $r_{Hn}$ is an $H$-radical property for $H$-module algebras.
  It is straightforward  to check  $r_{nH} \le r_{Hn}.$
 $\Box$

  $r_{Hn}$ is called the  $H$-von Neumann regular radical.

\begin {Theorem}\label {3.6}
If $I$ is an $H$-ideal of $R$, then $r_{Hn}(I) = r_{Hn}(R) \cap I $.
Namely, $r_{Hn}$  is a strongly hereditary $H$-radical property.

\end {Theorem}
{\bf Proof.}
By Lemma \ref {3.2}, $r_{Hn} (R) \cap I \subseteq r_{Hn}(I).$
Now, it is sufficient to show that $(x)_I = (x)_R$ for any $x \in r_{Hn}(I)$,
where $(x)_I$ and $(x)_R$ denote the $H$-ideals generated by $x$ in $I$ and
$R$, respectively. Let
$x = \sum (h_i \cdot x) b_i (h_i ' \cdot x)$ , where $h_i, h_i' \in H, b_i \in
I.$  We see that
\begin {eqnarray*}
R(H\cdot x) &=& R( H \cdot ( \sum (h_i \cdot x) b_i (h_i ' \cdot x)) \\
 &\subseteq &R (H\cdot x) I (H \cdot x)  \\
 &\subseteq & I (H \cdot x).
\end {eqnarray*}
Similarly, $$(H\cdot x)R \subseteq (H\cdot x)I.$$
Thus $(x)_I= (x)_R.$  $\Box$

A graded algebra $R$ of type $G$ is said to be Gr-regular if for every
homogeneous $a \in R_g$ there exists $b\in R$ such that $a= aba$
 \ \ \ ( see, \cite {NO} P258 ). Now, we give the relations between
 Gr-regularity and  $H$-regularity.

\begin {Theorem}\label {3.7}
              If $G$ is a finite group, $R$ is a graded algebra of type $G$,
               and $H=(kG)^*$, then    
$R$ is  Gr-regular iff  $R$ is  $H$-regular.

\end {Theorem}
{\bf Proof.} Let $\{ p_g \mid g\in G \}$ be the dual base of base $\{ g \mid g\in G
\}$. If $R$ is Gr-regular, for any
$a\in R$, then $a = \sum _{g \in G}  a_g$ with $a_g \in R_g$.
Since $R$ is Gr-regular, there exist $b_{g^{-1}} \in R_{g^{-1}}$
such that $a_g = a_g b_{g^{-1}} a_g $  and
$$ a = \sum _{g \in G} a_g = \sum _{g \in G } a_g b_{g^{-1}} a_g
= \sum _{g \in G}( p_g \cdot a ) b_{g^{-1}} (p_g \cdot a).$$
Consequently, $R$ is $H$-regular.

Conversely, if $R$ is $H$-regular then, for any $a \in R_g,$ there exist
$b_{x, y}\in R $  such that
$$a = \sum _{x,y \in G} (p_x \cdot a) c_{x,y} (p_y \cdot a).$$
Considering $a \in R_g$, we have that $a = a b_{g,g} a.$ Thus $R$ is Gr-regular
.$\Box$

\section { About J.R.Fisher's question }
 In this section,
 we answer the question  J.R.Fisher asked in \cite {F}. Namely,
 we give a necessary and sufficient condition for validity of relation (2) .

Throughout this section, let $k$ be a commutative ring with unit, $R$ an $H$-
module algebra and $H$ a Hopf algebra over $k.$

\begin {Theorem}\label {4.1}
Let  ${\cal K}$        be an ordinary special class of rings and closed with
respect to isomorphism. Set $r= r^{\cal K}$ and $\bar r_H (R) = r(R \# H) \cap R$
for any $H$-module algebra $R$. Then  $\bar r_H$ is an $H$-radical
property of
$H$-module algebras. Furthermore, it is an $H$-special radical.

\end {Theorem}

{\bf Proof.}
Let $\bar {\cal M}_R = \{ M \mid M $ is an $R$-prime module and $R/(0:M)_R \in
{\cal K } \}$    for any ring $R$ and $\bar {\cal M} = \cup \bar {\cal M}_R.$
Set ${\cal M}_R =\{ M \mid M \in \bar {\cal M }_{R \# H} \}$ for any $H$-
module algebra $R$ and ${\cal M}= \cup {\cal M}_R.$  It is straightforward to
check that $\bar {\cal M}$  satisfies the conditions of \cite [Proposition 4.3]
{Z97}. Thus ${\cal M}$ is an $H$-special module by      \cite [Proposition 4.3]
{Z97}.    It is clear that ${\cal M}(R) = \bar { \cal M} (R \#H) \cap R
= r(R \#H) \cap R$  for any $H$-module algebra $R$.
Thus $\bar r_H$ is an $H$-special radical by \cite [Theorem 3.1] {Z97}. $\Box$

Using the Theorem \ref {4.1}, we have that
all of $\bar r_{bH}, \bar r_{lH}, \bar r_{kH}, \bar r_{jH}, \bar r_{bmH}$
are $H$-special radical.

\begin {Proposition}\label {4.2}
Let ${\cal K}$ be a  special class of rings and closed with respect to
isomorphism.  Set $r= r^{\cal K}.$
Then

(1) $\bar r_H (R) \# H \subseteq r(R \#H);$

(2) $\bar r_H (R) \# H = r(R \#H)$  iff there exists an $H$-ideal $I$  of $R$
such that $r(R \#H) = I \#H$;

(3)  $R$ is an $ \bar r_H$-$H$-module algebra iff $r(R \#H) =R\# H$;

(4)  $I$ is an $ \bar r_H$-$H$-ideal of $R$ iff $r(I \#H)= I \#H$;

(5) $r(\bar r _H (R) \# H) = \bar r _H(R) \# H$ ;

(6) $r(R \# H) = \bar r _H(R) \# H$  iff
 $r(\bar r _H (R) \# H) =  r (R \# H).$

\end {Proposition}

{\bf Proof.} (1).  It is similar to the proof of Proposition \ref {2.3}.

(2). It is a straightforward verification.

(3). If $R$ is an $\bar r_H$-module algebra, then $R \# H \subseteq r(R \# H)$
 by part (1). Thus $R\# H = r(R \# H)$. The sufficiency is obvious.

(4), (5) and (6) immediately follow from part (3) . $\Box$

\begin {Theorem}\label {4.3}  If $R$ is an algebra over field $k$ with
unit and $H$ is a Hopf algebra over field $k$, then

(1)  $\bar r_{jH} (R) = r_{Hj}(R)$ and  $r_j(r_{Hj}(R)\# H) = r_{Hj}(R) \#H;$

(2)   $r_j(R\# H) = r_{Hj}(R) \#H$  iff
$r_j(r_{Hj}(R)\# H) = r_{j}(R \#H)$    iff
$r_j(r_j(R \# H)\cap R \# H) = r_j(R\# H);$

(3) Furthermore, if   $H$ is finite-dimensional, then
    $r_j(R\# H) = r_{Hj}(R) \#H$  { \ \ \ \ } iff \\
$r_j(r_{jH}(R)\# H) = r_{j}(R \#H).$

\end {Theorem}

{\bf Proof.}  (1) By \cite [Proposition 3.1 ]  {Z98}, we have
$\bar r_{jH} (R) = r_{Hj}(R)$.
Consequently, $r_j(r_{Hj}(R)  \# H) =r_{Hj}(R) \#H$ by Proposition
\ref {4.2} (5).

(2) It immediately  follows from part (1) and Proposition \ref {4.2} (6).

(3) It can easily be proved by part (2) and \cite [Proposition 3.3 (2)] {Z98}.
$\Box$

The theorem answers the question J.R.Fisher asked in \cite {F} :
When is $r_j(R \#H) = r_{Hj}(R)\#H$ ?

\begin {Proposition}\label {4.4}  If $R$ is an algebra over field $k$ with
unit and $H$ is a finite-dimensional Hopf algebra over field $k$, then

(1)  $\bar r_{bH} (R) = r_{Hb}(R)=r_{bH}(R)$
and  $r_b(r_{Hb}(R)\# H) = r_{Hb}(R) \#H;$

(2)   $r_b(R\# H) = r_{Hb}(R) \#H$  iff
$r_b(r_{Hb}(R)\# H) = r_{b}(R \#H)$  iff
$r_b(r_{bH}(R)\# H) = r_{b}(R \#H)$  iff
$r_b(r_{b}(R \#H) \cap R \# H) = r_{b}(R \#H).$

\end {Proposition}

{\bf Proof.}  (1). By \cite [Proposition 2.4 ]  {Z98}, we have
$\bar r_{bH} (R) = r_{Hb}(R)$.  Thus
$r_b(r_{Hb}(R)\# H) = r_{Hb}(R )\#H$  by Proposition \ref {4.2} (5).

(2). It follows from part (1) and Proposition \ref {4.2} (6) . $\Box$

In fact,  if $H$ is commutative or cocommutative, then $S^2 = id_H$
 by \cite [Proposition 4.0.1] {Sw},
and  $H$  is semisimple and cosemisimple iff
the character $char k $ of $k$ does  not
divides $dim H$   \ \ ( see, \cite [Proposition 2 (c)] {R94} ). It is clear that if $H$ is a finite-dimensional commutative or cocommutative Hopf algebra and
 the character $char k $ of $k$ does  not
divides $dim H$, then  $H$ is a finite-dimensional
semisimple and cosemisimple, commutative or cocommutative Hopf algebra.
Consequently, the conditions in Corollary \ref {2.9},
Theorem \ref {2.13}, Corollary  \ref{2.14} and \ref {2.17}
  can be simplified

\vskip 2cm

{\bf Acknowledgement }  I would like to express my gratitude
 to referee for his help.

\end {document}